\documentclass[a4paper,11pt]{article}
\usepackage[english]{babel}
\usepackage[psamsfonts]{amssymb}
\usepackage{cmmib57}
\usepackage{exscale}
\usepackage{amsmath}
\usepackage{amscd}
\usepackage{latexsym}
\usepackage{graphicx}
\usepackage{amsmath}
\usepackage{setspace}
\usepackage{fullpage}
\newtheorem{theorem}{Theorem}

\newtheorem{example}[theorem]{Example}

\newtheorem{lemma}[theorem]{Lemma}

\newtheorem{proposition}[theorem]{Proposition}
\newtheorem{remark}[theorem]{Remark}

\newenvironment{proof}[1][Proof]{\textbf{#1.} }{\ \rule{0.5em}{0.5em}}

\newcommand{\R}{\mathbb{R}}
\newcommand{\N}{\mathbb{N}}

\newcommand{\E}{\mathbb{E}}
\newcommand{\PP}{\mathbb{P}}

\newcommand{\un}{{\bf 1}}

\newcommand{\no}{\noindent}
\newcommand{\ds}{\displaystyle}
\newcommand{\vs}{\vspace}

\newcommand{\matrixnorm}[1]{\left|\!
\left|\!\left|#1\right|\!\right|\!\right|}
\def\darkbox{{\vrule height8pt width6pt depth0pt}}
\def\qed {\hfil \darkbox \break \vskip 0.25in}

\begin{document}
\baselineskip=1.1\baselineskip 
%
\title{On the geometric ergodicity of nonlinear
multivariate time series \vs{1cm}}
\date{}
\author{
Marco Ferrante
\footnote{corresponding author}
\\
Dipartimento di Matematica \\
Universit\`a degli Studi di Padova \\
via Trieste 63\\
35121 Padova, Italy \\
e-mail: ferrante@math.unipd.it
\and
Giovanni Fonseca
\\
Dip. di Scienze Economiche e Stat. \\
Universit\`a degli Studi di Udine
\\ via Tomadini, 30/A
\\ 33100 Udine, Italy \\
e-mail: giovanni.fonseca@uniud.it}
\maketitle
\begin{abstract}
\no
In this paper we consider multivariate time series
obtained as solution to multidimensional
nonlinear stochastic difference equations,
whose coefficients are
allowed to be locally degenerate and to present discontinuities.
We provide simple and easy to check
sufficient conditions for the
irreducibility, T-chain regularity and
geometric ergodicity of
these processes and apply the results to the BEKK-ARCH(1) models
with a nonlinear autoregressive term.

\end{abstract}

\vspace{0.2truecm}
{\bf Keywords:} Nonlinear stochastic difference equation,
T-chain, Geometric ergodicity, BEKK-ARCH(1).

\vspace{0.2truecm}
  {\bf AMS Classification: 60G10, 60J05, 62M10}

\vspace{0.2truecm}
{\bf Short title: Ergodicity of multivariate models}

\section{Introduction}

Let us consider a
system of nonlinear stochastic difference equations
\[
X_{t} = f(X_{t-1}) + g(X_{t-1}) e_t ,
\quad t\ge 1
\]
where $f:\R^n\rightarrow \R^n$,
$g:\R^n\rightarrow \R^{n\times k}$,
$\{e_t, t\in \N \}$ is a sequence
of independent, identically distributed $k$-dimensional
random vectors and $X_0$ is a given random vector.
We would like to find simple and easy to check conditions
that ensure the solution to be an irreducible, T-chain and
a geometric ergodic process.
The interest on these kind of models is clear:
first of all they can be thought as the discretization
of a multidimensional stochastic differential equation
and the properties of the discretized models are usually
of extreme interest.
Moreover, from a statistical point of view,
they can be considered as a state space representation
of a first order multivariate time series model and
many examples in the literature can be represented by these
models.

At a first sight this problem does not look 
very original and
worth to be studied in a new paper. However, under some very natural
conditions on the coefficients $f$ and $g$, like absence
of continuity and possible local singularity of
the matrix valued function $g$, to the best of our
knowledge no general results have been published so far,
except for those present in our former paper \cite{MR1977731},
in the case $n=k=1$.
It is worth to remark that several paper deal with
similar or more general models (see e.g. \cite{MR2341582},
\cite{MR1650873}, \cite{MR2488638}, \cite{MR2188304} and \cite{MR2352509}),
but in all
these papers stronger assumptions are required
on $g$, as everywhere continuity and non degeneracy.

The aim of this paper is to provide a first step in
order to fill this gap in the literature.
We will restrict ourselves to the case $n=k$
and we will assume that the noise random variables $e_t$
possess a strictly positive density on $\R^n$.
Under these conditions, we will be
able to obtain for this class of models
results similar to those that hold in the regular case,
adapting some of the standard techniques applied to
the smooth version of the present equation.
As an application, we shall consider a BEKK-ARCH(1)
multivariate model (see Engle and Kroner \cite{MR1325104},
Hansen and Rahbek \cite{hansen1998stationary} and Saikkonon
\cite{MR2352509})
and obtain a set of sufficient condition to be
this process geometric ergodic.

The paper is organized as follows:
in Section 2 we will present the model and
recall some notation and known results.
In Section 3 we will consider the problem to find out
sufficient conditions
to be the solution an aperiodic, irreducible, T-chain.
This part is fairly technical, but these three properties
are the fundamental ingredient in order to apply the
well known Foster-Lyapounov drift criteria of Section 4.
This technique allows us to determine a set of sufficient conditions
to be the solution an ergodic process.
In the last section we shall apply the results to a BEKK-ARCH(1) model
with a general nonlinear autoregressive term.

By $\lambda_n$ we will denote the Lebesgue measure on $\R^n$.
By $A^o$ we will denote the interior of the set $A$.
For $p\ge 1$, $\|\cdot\|_p$ will denote the $l_p$ norm on $\R^n$.
For $0<s\le 1$ and $x\in \R^n$, we will define
$\|x\|_s=\sum_{i=1}^n |x_i|^s$; this is clearly no more a norm,
but it still defines a pseudometric on $\R^n$, since
the triangular inequality holds true.
With $\matrixnorm{\cdot}$ we will denote
a generic matrix norm on $\R^{n\times n}$;
for $p\ge 1$,
$\matrixnorm{\cdot}_p$ will denote the
operator norm associated
with $\|\cdot\|_p$,
while
$\matrixnorm{\cdot}_{1,p}$ will denote the
{\em maximum column sum} matrix norm associated
with $\|\cdot\|_p$, whose
definition is as follows:
\[
\matrixnorm{A}_{1,p}=\max_{1\le j \le n} \|a_{\cdot j}\|_p
\quad.
\]
We will use the same definition when $0<s\le 1$, even if again
the function $\matrixnorm{\cdot}_{1,s}$ will not be a norm anymore.
Finally, by $\matrixnorm{\cdot}_{F}$ we will denote the Frobenius norm
\[
\matrixnorm{A}_F=\left(\sum_{1\le i,j \le n} |a_{i j}|^2\right)^{1/2}
\]
(see \cite{MR1084815}, Section 5.6 for a complete account
on this topic).

\no
For a map $G:X\times Y\rightarrow Z$, we shall denote by
$G_x$ the $x$-section of $G$, namely $G_x(y):=G(x,y)$, while,
given $B_1,\ldots,B_t\in \R^n$, we shall denote
$B_{1:t}=B_1\times \cdots \times B_t$
and, similarly, $u_{1:t}=(u_1,u_2,\ldots,u_t)\in \R^{nt}$.

\

\section{The multidimensional stochastic difference equation}

In this paper we will study
nonlinear stochastic difference equations
defined by the system
\begin{equation}
\left\{
\begin{array}{l}
\label{2.1}
X_{t} = f(X_{t-1}) + g(X_{t-1}) e_t ,
\quad t\ge 1
\\
X_0=\xi
\end{array}
\right.
\end{equation}
where $f:\R^n\rightarrow \R^n$,
$g:\R^n\rightarrow \R^{n\times n}$,
$\{e_t, t\in \N \}$ is a sequence
of independent, identically distributed $n$-dimensional
random vectors and $\xi$ is a given random vector.
In the discrete time case we have clearly no
problems related to the existence and uniqueness
of the solution, while a fundamental question is
if the system is ergodic, which is related to the fact that it
admits an invariant distribution.
We will be able to prove that, under some assumptions on
the coefficients $f$ and
$g$ and on the law of the noise sequence,
the solution
to (\ref{2.1}) is geometrically ergodic.

Let us start by considering the regularity of the coefficients $f$ and $g$.
In this paper, we consider the case with 
both the coefficients not everywhere
continuous and the matrix $g(x)$ 
locally singular.
This last assumption is not weird:
in the scalar case this means that the
function $g$ could be zero somewhere and this is indeed
the case if we choose $g$ to be an affine function.
Nevertheless in the literature
$g$ is usually assumed to be non singular and a possible
set of hypothesis
(see Liebscher \cite{MR2188304}) is that
there exist two constants $C_1,C_2>0$ such that
$\matrixnorm{g^{-1}(x)}\le C_1$ and
$|\det(g(x))|\le C_2$ on every compact
subset of $\R^n$ (see also Saikkonon
\cite{MR2352509} for the BEKK-ARCH model).

In the present paper, denoting by
$\Theta:=\{x\in \R^n: \det(g(x))\neq 0\}$
the set of ``regular'' points of $g$ and
by ${\cal C}_f$ (resp. ${\cal C}_g$) the set of
the continuity points of the function $f$
(resp. $g$),
we will require that the following two
assumptions are satisfied:
\begin{description}
\item{(H.1)}
$f$ and $g$ are
locally bounded
and the sets $\Theta$, ${\cal C}_f$ and ${\cal C}_g$ have not empty interior.
\end{description}

Under smooth conditions
on the coefficients $f$ and $g$, one fruitful approach is
to use the concept of the forward accessibility from the control
theory and its equivalence with the much more workable Rank condition
(see e.g. Meyn-Tweedie \cite{MR2509253}, Chapter 7).
Since in our case we do not assume differentiability
of the coefficients, we have to find out
a different approach, even if our property will
be at the end stronger, but not so far from the
forward accessibility.
Let us denote by $O\subseteq \R^n$ the support of
the random vector $e_1$,
$F(x,u)=f(x)+g(x)u$ and, inductively, for $t\in \N^+$
\begin{equation}
\label{iter}
F^{t+1}(x_0,u_1,\ldots,u_{t+1}):=
F(F^t(x_0,u_1,\ldots,u_t),u_{t+1})
\quad.
\end{equation}
When $t=0$, $F^t(x_0,u_1,\ldots,u_t)\equiv x_0$.

Our second assumption will be that
\begin{description}
\item{(H.2)}
For any $x_0\not\in \Theta$, there exists $t\in \N^+$ and
$u_1,\ldots,u_{t}\in O
$
such that
$F^{t}(x_0,u_1,\ldots,u_t) \in (\Theta \cap {\cal C}_f \cap {\cal C}_g)^o$ and
$F^{t-1}(x_0,u_1,\ldots,u_{t-1}) \in {\cal C}_f \cap {\cal C}_g$.
\end{description}
\begin{remark}
\label{fa}
Under (H.2) and assuming that $O$ has non-empty interior,
we easily obtain that
for any $x_0$, there exists $t\in \N^+$ and
$
u_1,\ldots,u_{t}\in O
$
such that
$F^{t+1}_{x_0,u_1,\ldots,u_t}(O)$ has non-empty interior.
This condition implies, but is evidently
stronger than the {\em forward accessibility}, which
requires that for any $x_0\in \R^n$,
$\bigcup_{t=0}^{+\infty} F^t_{x_0}(O^t)$
has non-empty interior.
\end{remark}

\begin{remark}
\label{sc}
If we assume that $f$ and $g$ are continuous, a sufficient condition
in order to
satisfy assumption (H.2) is that for any $x\in \R^n$, there exists
$t\in \N$ such that
$ f^t(x) \in \Theta^o$.
\end{remark}

To conclude, let us state the assumptions on the noise sequence
$\{e_t, t\in \N \}$.
The price to be paid for the
weak hypothesis (H.1) is quite expensive, since we
have to assume absolute continuity and
lower semicontinuity of the noise density.
However, this is often the additional condition that we have to
ask in order to allow some kind of singularity in the
coefficient $g$ (see e.g the results on the bilinear processes
in \cite{MR648732} and \cite{MR1892989}, the threshold bilinear processes
in \cite{MR1977731}
and the nonlinear state space models in \cite{MR2509253}).
\begin{description}
\item{(H.3)}
$\{e_t, t\in \N \}$ is a sequence
of independent, identically distributed $n$-dimensional random vector,
absolutely continuous w.r.t. Lebesgue measure $\lambda_n$ on
${\cal B}(\R^n)$, with
density $p(\cdot)$ strictly positive almost everywhere and lower
semicontinuous.
\end{description}

\section{Irreducibility, aperiodicity and T-chain property}

In order to apply the classic Foster-Lyapunov drift criteria
of the next section (see
Meyn-Tweedie \cite{MR2509253} for a comprehensive introduction
to this topic), we need three basic ingredients.
Indeed, we have to prove that the Markov chain, solution to
(\ref{2.1}),
is $\varphi$-irreducible, for a given measure $\varphi$, aperiodic and
a T-chain.
While it is usually not too difficult to find out a set of reasonable
conditions that ensure that the process is $\varphi$-irreducible
and aperiodic, it is more challenging to handle the T-chain
condition.
Most of the papers in the literature do not spend much time on this
part of the study and the authors usually state some general
conditions that ensure the process to satisfy these three properties
(see e.g. Liebscher \cite{MR2188304}, Section 4).
Since in this paper we would like to allow the diffusion coefficient $g$
to be locally singular, we shall need to impose the previous set of
stronger assumptions (H.1)-(H.3).

Let us start by stating a simple result, that we will need
in the sequel and whose proof is immediate (see \cite{MR1977731}).
\begin{lemma}
\label{lsc}
Let $A,B\subseteq \R^n$:
if $F:A \rightarrow B$ is a continuous function
and $G,H:B \rightarrow \R$ are two lower semicontinuous (lsc)
functions, then
$G\circ F$ and $G \cdot H$ are lsc.
\end{lemma}

We are now able to prove the main result of this section:
\begin{proposition}
\label{propTchain}
Under $(H.1)-(H.3)$, the process solution to
(\ref{2.1}) is a $\lambda_n$-irreducible, aperiodic T-chain.
\end{proposition}

\proof {\bf $\lambda_n$-irreducibility}: we have to prove that for any $A\in
{\cal B}(\R^n)$, such that $\lambda_n(A)>0$, and any $x\in \R^n$, there
exists $t=t(x,A)\in \N$ such that $P^t(x,A)=\PP [X_t\in A|X_0=x]>0$.
If $\det(g(x))\neq 0$, we get $P(x,A)>0$.
Otherwise, by assumption
(H.2) we get that there exists $t\in \N$ and $u_1,\ldots,u_t\in \R^n$
such that $F^t(x,u_1,\ldots,u_t)$,
defined in (\ref{iter}), belongs to $(\Theta \cap {\cal C}_g)^o$, and is
continuous in $(x,u_1,\ldots,u_t)$.
Therefore, there exist open balls $B_1,\ldots,B_t$ in
$\R^n$ such that
$F^t(x,B_1,\ldots,B_t)\subseteq \Theta^o$.
By
(H.1)-(H.3) we get
\[
\begin{array}{l}
{\ds P^{t+1}(x,A) =
\PP\left[f(F^t(x,e_{1:t}))+g(F^t(x,e_{1:t}))e_{t+1}\in A\right]}
\\ \\
\ \ \ \ \
{\ds \ge \int_{B_{1:t}}
\Big[
\int_A
\left|\det (g(F^t(x,u_{1:t})))\right|^{-1}
p(g(F^t(x,u_{1:t}))^{-1}(u_{t+1}
-f(F^t(x,u_{1:t})))) du_{t+1}
\Big] \times}
\\ \\
{\ds
\ \ \ \ \
p(u_1)\cdots p(u_t)
d u_{1:t}
\ge
c_3 \int_A \Big[ \int_{B_{1:t}}
p(g(F^t(x,u_{1:t}))^{-1}(u_{t+1}
-f(F^t(x,u_{1:t}))))
\times}
\\ \\
{\ds
\ \ \ \ \
p(u_1)\cdots p(u_t)
d u_{1:t}
\Big]  du_{t+1} > 0 \quad ,}
\end{array}
\]
where $c_3:=\inf_{u_{1:t}\in B_{1:t}}\left|
\det (g(F^t(x,u_{1:t})))\right|^{-1}< \infty$,
and the $\lambda_n$-irreducibility is proved.

{\bf Aperiodicity}: we will prove that the solution process is strongly
aperiodic, i.e. that there exist a nontrivial measure $\nu_1$ on
${\cal B}(\R^n)$ and a subset $A\in {\cal B}(\R^n)$, with
$\nu_1(A)>0$, such that for any $x \in A$ and $B\in {\cal
B}(\R^n)$, $P(x,B)\ge \nu_1(B)$.
Let us take $x\in (\Theta \cap {\cal C}_f\cap {\cal C}_g)^o$; by the
assumption (H.1) we get that there exists an open
bounded neighborhood $A$ of
$x$ and two positive constants $c_1, c_2$ such that
\[
0 < c_1 \le \left|\det(g(y))\right| \le c_2
\quad
\]
for any $y\in A$.
By (H.2) we get

\[
\begin{array}{l}
{\ds P(x,B) = \int_B \left|\det (g(x))\right|^{-1} p(g(x)^{-1}(y
-f(x))) dy}
\\ \\
\ \ \ \ \ \ \ \ \ \ \
{\ds \ge c_2^{-1} \int_{A\cap B} p(g(x)^{-1}(y -f(x))) dy
\ge c_2^{-1} k_1 \lambda_n(A\cap B)} \quad,
\end{array}
\]
where ${\ds 0<k_1=\inf_{x,y\in A} p(g(x)^{-1}(y -f(x)))}$. The result
holds for $\nu_1(\cdot)= c_2^{-1} k_1 \lambda_n(A\cap
\cdot)$.

{\bf T-chain condition}: By Proposition 6.4.2 in Meyn and Tweedie
\cite{MR2509253}, it will be sufficient to prove that for each $x \in
\R^n$, there exists a $t \in \N$ and a non trivial substochastic
transition kernel $T_x(\cdot,\cdot)$. l.s.c. in the first variable,
such that $P^t(y,A)\ge T_x(y,A)$ for each $y\in \R^n$ and $A\in {\cal
B}(\R^n)$.
Let $x\in \R^n$: by assumption
(H.2) we get that there exist $t\in \N$ and $u_1,\ldots,u_t\in \R^n$
such that $F^t(x,u_1,\ldots,u_t)\in (\Theta \cap {\cal C}_f\cap {\cal C}_g)^o$ and
$t+1$ open bounded sets $B_0,B_1,\ldots,B_t$ of $\R^n$,
such that
$F^t(B_0,B_1,\ldots,B_t)\subseteq \Theta^o$.
Moreover, we can assume that
$f$ and $g$ are continuous on
$F^t(B_0,B_1,\ldots,B_t)$.
Hence, for $y\in B_0$ and $A\in {\cal B}(\R^n)$, we get
\[
\begin{array}{l}
{\ds
P^{t+1}(y,A) =
\PP\left[f(F^t(y,e_{1:t}))+g(F^t(y,e_{1:t}))e_{t+1}\in A\right]}
\\ \\
\ \ \ \ \
\ge
{\ds \int_{B_{1:t}}
\left[
\int_A
\left|\det (g(F^t(y,u_{1:t})))\right|^{-1}
p(g(F^t(y,u_{1:t}))^{-1}(u_{t+1}
-f(F^t(y,u_{1:t})))) du_{t+1}
\right] \times}
\\ \\
\ \ \ \ \
{\ds
p(u_1)\cdots p(u_t)
d u_{1:t}
\ge
c_4 \int_A \left[ \int_{B_{1:t}}
p(g(F^t(y,u_{1:t}))^{-1}(u_{t+1}
-f(F^t(y,u_{1:t}))))
\right. \times}
\\ \\
\ \ \ \ \
{\ds
p(u_1)\cdots p(u_t)
d u_{1:t}
\Big]  du_{t+1}
=: \widetilde{T}(y,A)}
\quad ,
\end{array}
\]
where $c_4:=\inf_{(y,u_{1:t})\in B_{0:t}}\left|
\det (g(F^t(y,u_{1:t})))\right|^{-1}<\infty$.
By Lemma \ref{lsc} and Fatou's Lemma we get that
$\widetilde{T}(y,A)$ (for $y \in B_0$) is a lsc function and we
can define the substochastic
transition kernel
$T_x(y,A):=\phi(y) \widetilde{T}(y,A)$,
with $\phi(\cdot)$ a smooth function
whose support is contained in $B_0$.
It is clear that
$P^t(y,A)\ge T_x(y,A)$ for each $y\in \R^n$ and
$A\in {\cal B}(\R^n)$ and
the proof is complete.
\qed
We conclude this section by considering a simple bivariate
time-series,
solution of a two dimensional difference equation,
where is present a threshold and a singular part.

\begin{example}
\label{model-expol_1}
Let us take $n=2$ and consider the difference equation
(\ref{2.1}), with $f$ and $g$ defined as follows
\[
f(x,y)=
\left(
  \begin{matrix}
  a_1 \\
  a_2
  \end{matrix}
\right)
+
\left(
  \begin{matrix}
  b_{11} & b_{12} \\
  b_{21} & b_{22}
  \end{matrix}
  \right)
\left(
  \begin{matrix}
  x \\
  y
  \end{matrix}
\right)
\]
and
\[
g(x,y)=
\left(
  \begin{matrix}
  d_{11} \cdot x & d_{12} \cdot y \\
  d_{21} \cdot x & d_{22} \cdot y
  \end{matrix}
  \right)
  \un_{\R^2 \setminus C}
+
\left(
  \begin{matrix}
  d_{31} \cdot x & 0 \\
  d_{32} \cdot y & 0
  \end{matrix}
  \right)
  \un_{C}
+
\left(
  \begin{matrix}
  d_{41} & 0 \\
  d_{42} & 0
  \end{matrix}
  \right)
\]
where $C=\{x\le 0,y\le 0\}$ and $d_{11} d_{22}-d_{12} d_{21}\neq 0$.
Note that the matrix $g(x,y)$ is singular for $(x,y)\in C\cup D_1$ and
could be singular for $(x,y)\in D_2$, where $D_1:=\{(x,y):x=0, y>0\}$
and $D_2:=\{(x,y):x>0, y=0\}$. Moreover, $g$
is not continuous on the boundary of $C$.
It is easy to prove that if $f$ and $g$ are such that
\begin{equation}
\label{coefexpol}
d_{11} d_{42}-d_{21} d_{41}\neq 0
\ , \
d_{31} d_{42}- d_{32} d_{41}\neq 0
\end{equation}
and the noise sequence admits $\R^n$ as its support, then
the hypotheses (H.1)-(H.2) are satisfied.
Indeed, the only
non trivial part is (H.2), which is always satisfied
since,
for $(x,y)\in C$,
there exists
$u\in \R^2$ such that
$f(x,y)+g(x,y) u\not\in C\cup D_1\cup D_2$,
while
for $(x,y)\in D_i$, $i=1,2$,
there exist
$u, v\in \R^2$ such that
$f(f(x,y)+g(x,y) u)+g((f(x,y)+g(x,y) u))v\not\in D_i$.

\end{example}

\section{Geometric ergodicity}

In this section we will obtain, in a standard way, a set of sufficient
conditions on the coefficients of the difference equation in order to be
the solution process geometrically ergodic.
Due to the weak assumptions on this coefficients, we will obtain
a rather strict sufficient condition, but this is in line with
previous results in the literature. Since our approach to prove the
geometric ergodicity is based
on the choice of a drift function, for specific models
like the BEKK-ARCH(1) model in the next section, it
could be more convenient to use a different function, but
the proof will be similar to the one presented 
here.

Let us consider the process $\{X_t, t\ge 0\}$ solution to
equation (\ref{2.1}) and assume that it is a
$\lambda_n$-irreducible, aperiodic T-chain.
In order to apply the classic Foster-Lyapunov drift criteria
for $V(x)=1+\|x\|_s$, when $s>0$,
we will need to apply some easy properties of the
functions $\|\cdot\|_s$ and $\matrixnorm{\cdot}_{1,s}$.
Given $A,B$ two $n \times n$ real matrices and
$x\in \R^n$,
it holds that
\[
\matrixnorm{A B}_{1,s} \le \matrixnorm{A}_{1,s} \matrixnorm{B}_{1,s}
\]
and that
\[
\|Ax\|_s\leq \matrixnorm{A}_{1,s} \|x\|_s
\quad,
\]
for any $s>0$; these results are well known for $s\ge 1$ and
immediate to be proven for $s<1$.
Note that, the same properties hold true for
$\matrixnorm{\cdot}_s$ instead of
$\matrixnorm{\cdot}_{1,s}$, when $s\ge 1$, and
furthermore, when $s=2$, for $\matrixnorm{\cdot}_F$ instead of
$\matrixnorm{\cdot}_{1,2}$.

We are now ready to provide an easy to check
set of sufficient conditions in order to be the solution
process geometric ergodic.

\begin{proposition}
\label{propn1}
Let $ \{ X_t \} $ be
the solution process of (\ref{2.1}) and
assume that it is $\lambda_n$-irreducible, aperiodic, T-chain.
If for $s>0$
\begin{enumerate}
\item
$f$ and $g$ are locally bounded;
\item
there exist $a_f \geq 0$ and $ M, a_g, b_f, b_g > 0$ such that
\[
\|f(x)\|_s \le a_f+b_f \|x\|_s
\ \ , \ \
\matrixnorm{g(x)}_{1,s} \le a_g+b_g\|x\|_s
\]
for any $x \in \R^n$ with $\|x\|_s>M$;
\item
$ \gamma = b_f \ + \ b_g\ \E\left[\|e_1\|_s\right] <1 \ ;$
\end{enumerate}
then $ \{ X_t, t\in \N \} $ is geometrically ergodic.

\no
Furthermore, if the previous conditions hold for $s\ge 1$, then
each component of the stationary distribution has finite moments up to order $s$.

\end{proposition}


\proof
Let us consider the function $V(x)=1+\|x\|_s$ for an arbitrary $s>0$
and the compact set $C=\{x \in \R^n : \|x\|_s\le M\}$.
Since the solution is a $\lambda_n$-irreducible T-chain,
we have that $C$ is also {\em petite}.

By triangular inequality and assumptions i.--iii. we obtain
\begin{eqnarray}
\E[V(X_t)|X_{t-1}=x] 
& \leq & (b_f + b_g \E[\|e_t\|_s])(1+ \|x\|_s) +\nonumber \\
& + & a_f+a_g \E[\|e_t\|_s]+1-(b_f + b_g \E[\|e_t\|_s])\nonumber
\end{eqnarray}
for every $x \in \R^n$ such that $\|x\|_s>M,$ and
\begin{eqnarray}
\E[V(X_t)|X_{t-1}=x] & \leq & b_{M}<\infty \quad \forall \ x \in C.
\nonumber
\end{eqnarray}
Summarizing, for any $x \in \R^n$ it holds
\begin{eqnarray}
\E[V(X_t)|X_{t-1}=x] & \leq & (b_f + b_g \E[\|e_t\|_s]) V(x) + b_M \un_{C}.
\nonumber
\end{eqnarray}

If $b_f + b_g \E[\|e_t\|_s]<1$, applying
Lemma 15.2.8,
Theorem 15.0.1
and Theorem 14.0.1
in Meyn and Tweedie \cite{MR2509253}, we
get that $\{X_t\}$ is  a
geometrically ergodic Markov chain
and, when $s\ge 1$, that the moments
of the components of the stationary distribution
are finite up to order $s$. \qed

\begin{remark}
\label{fro}
A simple extension to the previous result can be obtained
by taking, for $s\ge 1$, $\matrixnorm{\cdot}_s$ or $\matrixnorm{\cdot}_F$
instead of $\matrixnorm{\cdot}_{1,s}$.
More generally, any norm on $\R^n$
(and the corresponding induced matrix norm) could be considered.
Indeed, in the proof the only property
of the matrix norms that we use is 
$\|g(x)e_t\|_s \leq \matrixnorm{g(x)}_{1,s}\|e_{t}\|_s$ which holds true for these matrix norms too.
Our choice of the matrix norm $\matrixnorm{\cdot}_{1,s}$
is justified by the possibility to use in assumption iii. absolute moments of order smaller then 1, which weakens
the restriction on the noise.
\end{remark}

\begin{example}
\label{model-expol_2}
Let us consider the threshold model of Example \ref{model-expol_1}.
We will assume that the i.i.d. random sequence
$\{e_t,t\in \N\}$ is distributed
according to an Expol(2) law (see \cite{MR1381677})
with density
$ p(x,y) \propto \exp(-(x^2-1)^2 -(y^2-1)^2)$.
This density
presents four global maxima
and one local minimum and is straightforward to simulate by
standard methods. 
In order to apply the results in Proposition \ref{propn1},
it is easy to check that Assumption ii. hold true,
when $s=1$, for
\[
b_f=\max\{|b_{11}| + |b_{21}|,|b_{12}| + |b_{22}|\}
\]
and
\[
b_g=\max\{|d_{11}| + |d_{21}|, |d_{31}|,|d_{32}|,|d_{12}| + |d_{22}|\}
\]
Since $\E[\|e_1\|_1]\sim 1.66$, it easy to see that the previous
model with
\[
f(x,y)=
\left(
  \begin{matrix}
  a_1 \\
  a_2
  \end{matrix}
\right)
+
\left(
  \begin{matrix}
  0.2 & 0.1 \\
  0.1 & 0.3
  \end{matrix}
  \right)
\left(
  \begin{matrix}
  x \\
  y
  \end{matrix}
\right)
\]
and
\[
g(x,y)=
\left(
  \begin{matrix}
  0.1 x &  -0.15 y \\
  -0.15 x & 0.1 y
  \end{matrix}
  \right)
  \un_{\R^2 \setminus C}
+
\left(
  \begin{matrix}
  0.2 x & 0 \\
  -0.25 y & 0
  \end{matrix}
  \right)
  \un_{C}
+
\left(
  \begin{matrix}
  1 & 0 \\
  1 & 0
  \end{matrix}
  \right)
\]
satisfies the assumptions of Proposition \ref{propn1}, with
$b_f \ + \ b_g\ \E\left[\|e_1\|_1\right]\sim 0.981$.
On the other hand, if we modify the coefficients for the conditional mean,
for example taking $b_{11}=b_{22}=1, b_{12}=b_{21}=0 $,
the previous assumptions are no more satisfied and the
distribution of the solution process does not converge to any distribution,
as can be seen by simulation.
The same result is obtained if we modify the coefficients for the conditional variance,
for example taking $d_{11}=d_{12}=d_{21}=d_{22}=0.4$.
Nonetheless simulation of the limit distribution for
other set of values shows that even if the
sufficient condition is (slightly) violated, the model
still appear ergodic, but as pointed out before, Assumption iii.
is strong.
\end{example}

\section{Multivariate BEKK-ARCH(1) models with nonlinear autoregressive terms}

In this final section we will consider a
locally degenerate multivariate BEKK-ARCH(1) model,
with a nonlinear autoregressive term.
This model belongs to the multivariate BEKK-GARCH class,
first proposed by Engle and Kroner in \cite{MR1325104},
which is particularly useful in multivariate financial time-series,
since allow to model both the variances and the covariances.
Contrary to all previous works, we will ask that
the matrix valued coefficient will be just positive
semidefinite and we will be able
to derive simple sufficient conditions for the
regularity and geometric ergodicity of the solution process.

Let us consider a process
$ \{ X_t, t\ge 0 \} $, solution to the following
difference equation:
\begin{equation}
\label{BEKK}
X_{t} = f(X_{t-1}) + (B+(AX_{t-1}) (AX_{t-1})^T)^{1/2} e_t ,
\quad t\ge 1
\end{equation}
where $f:\R^n\rightarrow \R^n$,
$B$ is a positive semidefinite $n\times n$ real matrix,
$A$ is a $n\times n$ real matrix and
$\{e_t, t\in \N \}$ is a sequence
of independent, identically distributed $n$-dimensional
random vectors.
Since $(Ax) (Ax)^T$ is, for any $x\in \R^n$, a
positive semidefinite $n\times n$ real matrix, $B+(Ax) (Ax)^T$ is itself a positive semidefinite matrix and there exists a unique, positive
semidefinite square root
(see \cite{MR1084815}, Chapter 7).
The process $X_t$, solution to (\ref{BEKK}) is called in the
time-series literature a BEKK-ARCH(1) model and its ergodicity
has been considered just for the regular case, i.e.
assuming that the matrix $B+(Ax) (Ax)^T$ is positive definite
and its smaller eigenvalue is uniformly bounded away from zero
(see Saikkonon \cite{MR2352509}).

To determine a set of sufficient conditions to be the solution process
a $\lambda_n$-irreducible, aperiodic, T-chain, we will
consider for simplicity the case $n=2$.
Let us assume that the matrix $B$ will be non zero and
positive semidefinite, but not
positive definite.
In this case the matrix $(B+(Ax)(Ax)^T)^{1/2}$ will be
positive semidefinite, but could be not positive definite.
Indeed, the determinant
of the matrix valued function $g(x)=B+(Ax)(Ax)^T$ is equal to
\[
\det(g(x))=b_{11} (a_{21}x_1 + a_{22}x_2)^2 + b_{22} (a_{11}x_1 + a_{12}x_2)^2
- 2 b_{12} (a_{21}x_1 + a_{22}x_2)(a_{11}x_1 + a_{12}x_2)
\]
and a simple computation show that
this determinant is zero for
any $(x_1,x_2)\in \R^2$ if
both $a_{11} b_{22}^{1/2} - a_{21} b_{11}^{1/2}$ and
$a_{12} b_{22}^{1/2} - a_{22} b_{11}^{1/2}$ are zero, while
otherwise
it is zero just on a straight line $L=\{(x_1,x_2)\in \R^2, \alpha x_1 = \beta x_2\}$,
for suitable constants $\alpha$ and $\beta$.
With the notation of Section 2, we get that in the latter case the set of
regular points of $g$, $\Theta$, coincides with $L^c$.
Moreover, the function
$x\rightarrow(B+(Ax)(Ax)^T)^{1/2}$ is continuous if
$B+(Ax)(Ax)^T$ is positive definite.

A set of sufficient conditions for the assumptions (H.1) and (H.2) will
be as follows:
\begin{description}
\item{(B.1)}
At least one between $a_{11} b_{22}^{1/2} - a_{21} b_{11}^{1/2}$ and
$a_{12} b_{22}^{1/2} - a_{22} b_{11}^{1/2}$ is different of zero and
for any $x\in L=\{x\in \R^2: \det(g(x))=0\}$, there exist $u,v\in \R$
such that $y=f(x)+g(x)u \in {\cal C}_f\cap L^c$ and
$f(y)+g(y)v \in ({\cal C}_f\cap L^c)^o$, with
${\cal C}_f^o \neq \emptyset$.
\end{description}
\begin{remark}
Clearly, for a specific choice of the autoregressive term $f$,
one can provide better sufficient conditions in
order to be the assumptions (H.1) and (H.2) satisfied.
\end{remark}
The next result follows as a simple corollary
of previous Proposition \ref{propTchain}:
\begin{proposition}
\label{propTchainBEKK}
Let $n=2$ and assume that $(B.1)$ and $(H.3)$ are satisfied.
Then, the process solution to
(\ref{BEKK}) is a $\lambda_n$-irreducible, aperiodic T-chain.
\end{proposition}

Let us now consider the geometric ergodicity:
in order to apply the results of the previous sections,
we will use here the Frobenius matrix norm. In fact we will make use
of the fact that
$\matrixnorm{A}_F=\left(\sum_{i,j=1}^n a_{ij}^2\right)^{1/2}
= \left(tr (A^T A)\right)^{1/2}$, which gives in our case that
\[
\matrixnorm{(B+(Ax) (Ax)^T)^{1/2}}_F^2=
tr (B) + tr((Ax) (Ax)^T) .
\]
A set of sufficient condition for the
geometrically ergodicity of the present model
follows as a simple modification of previous Proposition 6,
whose simple proof we omit.
\begin{proposition}
\label{propn1BEKK}
Let $ \{ X_t \} $ be
the solution process of (\ref{BEKK}) and
assume that it is a $\lambda_n$-irreducible, aperiodic, T-chain.
If
\begin{enumerate}
\item
$f$ is locally bounded;
\item
there exist $a_f \geq 0$ and $ M, b_f > 0$ such that
\[
\|f(x)\|_2 \le a_f+b_f \|x\|_2
\]
for any $x \in \R^n$ with $\|x\|_2>M$;
\item
$ \gamma = b_f \ + \ \matrixnorm{A}_F \ \E\left[\|e_1\|_2\right] <1 \ ;$
\end{enumerate}
then $ \{ X_t, t\in \N \} $ is geometrically ergodic.

\no
Furthermore,
the components of the stationary distribution have finite second moments.

\end{proposition}

\begin{remark}
\label{rho}
More general conditions for geometric ergodicity are present in the literature of the BEKK-ARCH models (see, for instance, \cite{MR2352509}), but a basic ingredient
of all the proofs is that the matrix $B+(Ax) (Ax)^T$ will be positive definite.
\end{remark}

\bibliography{HT}  
\bibliographystyle{plain}

\end{document}